\documentclass[10pt]{amsart}
\usepackage[cp1251]{inputenc}
\usepackage[english,russian]{babel}
\usepackage{amsmath}
\usepackage{longtable}
\usepackage{tikz}
\usepackage{amssymb}
\usepackage{amsfonts}

\newtheorem{lemma}{Lemma}
\newtheorem{thm}{Theorem}

\newtheorem{prop}{Proposition}

\newtheorem{ex}{Example}

\def\top
     {

     \begin{tabular}{c}
     \hphantom{aaaaaaaaaaaaaaaaaaaaaaaaaaaaaaaaaaaaaaaaaaaaaaaaaaaaaaaaaaaaaaaaaaaaaa} \\
     \end{tabular}

   {\flushleft
   \hspace{65mm}{\rm\small 
   } 
     \newline
         \hphantom{aaaaaaaaaaaaaaaaaaaaaaaaaaaaaaaaaaaaaaaaaaaaaaaaaaaaaaaa}{\rm\small MSC\ \ 20D60,~20D06 }
  }
}

\begin{document}
\renewcommand{\refname}{References}
\renewcommand{\proofname}{Proof.}
\thispagestyle{empty}

\title[On the Pronormality of Subgroups of Odd Index]{ON THE PRONORMALITY OF SUBGROUPS OF ODD INDEX IN FINITE SIMPLE GROUPS}
\author{{Anatoly~S.~Kondrat'ev}}%
\address{Anatoly~Semenovich~Kondrat'ev
\newline\hphantom{iii} Krasovskii Institute of Mathematics and Mechanics of the Ural Branch of the  Russian Academy of Science,
\newline\hphantom{iii} 16, S. Kovalevskaja Street,
\newline\hphantom{iii} 620990, Yekaterinburg, Russia
\newline\hphantom{iii} Ural Federal University named after the first President of Russia B.~N.~ Yeltsin,
\newline\hphantom{iii} 19, Mira Street,
\newline\hphantom{iii} 620002, Yekaterinburg, Russia}
\email{a.s.kondratiev@imm.uran.ru}
\author{{Natalia~V.~Maslova}}%
\address{Natalia~Vladimirovna~Maslova
\newline\hphantom{iii} Krasovskii Institute of Mathematics and Mechanics of the Ural Branch of the  Russian Academy of Science,
\newline\hphantom{iii} 16, S. Kovalevskaja Street,
\newline\hphantom{iii} 620990, Yekaterinburg, Russia
\newline\hphantom{iii} Ural Federal University named after the first President of Russia B.~N.~ Yeltsin,
\newline\hphantom{iii} 19, Mira Street,
\newline\hphantom{iii} 620002, Yekaterinburg, Russia}
\email{butterson@mail.ru}
\author{{Danila~O.~Revin}}%
\address{Danila~Olegovich~Revin
\newline\hphantom{iii} Sobolev Institute of Mathematics of the Siberian Branch of the  Russian Academy of Science,
\newline\hphantom{iii} 4, Koptuga Avenue,
\newline\hphantom{iii} 630090, Novosibirsk, Russia
\newline\hphantom{iii} Novosibirsk State University,
\newline\hphantom{iii} 1, Pirogova Street,
\newline\hphantom{iii} 630090, Novosibirsk, Russia}
\email{revin@math.nsc.ru}


\top \vspace{1cm}
\maketitle {\small
\begin{quote}
\noindent{\sc Abstract. } A subgroup $H$ of a group $G$ is said to be {pronormal} in $G$ if $H$ and $H^g$ are conjugate in $\langle H, H^g \rangle$ for every $g \in G$. Some problems in finite group theory, combinatorics, and permutation group theory were solved in terms of pronormality. In 2012, E.~Vdovin and the third author conjectured that the subgroups of odd index are pronormal in finite simple groups.
In this paper we disprove their conjecture and discuss a recent progress in the classification of finite simple groups in which the subgroups of odd index are pronormal.
\\

 \end{quote}
}

Throughout the paper we consider only finite groups, and thereby the term ''group'' means ''finite group''.

According to P. Hall, a subgroup $H$ of a group $G$ is said to be {\it pronormal} in $G$ if $H$ and $H^g$ are conjugate in $\langle H, H^g \rangle$ for every $g \in G$.

Some of well-known examples of pronormal subgroups are the following: normal subgroups; maximal subgroups; Sylow subgroups; Sylow subgroups of proper normal subgroups; Hall subgroups of solvable groups. In 2012, E.~Vdovin and the third author \cite{Vdovin_Revin} proved that the Hall subgroups are pronormal in all simple groups.

The following assertion give a connection between pronormality of subgroups and properties of permutation representations of finite groups.
\begin{thm}{\rm{\bf (Ph. Hall, 1960s).}}\label{Halls1Thm} Let $G$ be a group and $H \le G$.
$H$ is pronormal in $G$ if and only if in any transitive permutation representation of $G$, the subgroup $N_G(H)$ acts transitively on the set $fix(H)$ of fixed points of $H$.
\end{thm}

Pronormality is the universal property with respect to Frattini Argument. Indeed, it is not hard to prove the following proposition (see \cite[Lemma~4]{Guo_Mas_Rev}).

\begin{prop}\label{Frattini} Let $G$ be a group, $A\unlhd G$, and $H\leq A$. Then the following statements are equivalent{\rm:}

$(1)$ $H$ is pronormal in $G${\rm;}

$(2)$ $H$ is pronormal in $A$ and $G=AN_G(H)${\rm;}

$(3)$ $H$ is pronormal in $A$ and $H^A=H^G$.

\end{prop}

\smallskip

Some problems in finite group theory, combinatorics, and permutation group theory were solved in terms of pronormality.
 For example, according to L.~Babai \cite{Babai}, a group $G$ is called a {\it CI-group} if between every two isomorphic relational structures on $G$ (as underlying set) which are invariant under the group $G_R=\{g_R\mid g\in G\}$ of  right multiplications
 $g_R:x\mapsto xg$ (where $g, x \in G$), there exists an isomorphism which is at the same time an automorphism of~$G$. Babai \cite{Babai} proved that a group $G$ is a CI-group if and only if $G_{R}$ is pronormal in $Sym(G)$. In particular, if $G$ is a CI-group, then $G$ is abelian. With using mentioned Babai's result, P.~Palfy \cite{Palfy} obtained a classification of CI-groups.

 \medskip

Thus, the following problem naturally arises.

\medskip

{\bf General Problem.} {\it Given a finite group $G$ and $H \le G$, is $H$ pronormal in $G$?}

\medskip

\medskip

Ch.~Praeger \cite{Praeger} investigated pronormal subgroups of permutation groups. She proved the following theorem.

\begin{thm}\label{PrThm}  Let $G$ be a transitive permutation group on a set $\Omega$ of $n$ points, and let $K$ be a non-trivial pronormal subgroup of $G$. Suppose that $K$ fixes exactly $f$ points of $\Omega$. Then $f \le \frac{1}{2}(n-1)$, and if $f = \frac{1}{2}(n-1)$, then $K$ is transitive on its support in $\Omega$, and either $G \ge Alt(n)$, or
$G = GL_d(2)$ acting on the $n = 2^d - 1$ non-zero vectors, and $K$ is the pointwise stabilizer of a hyperplane.
\end{thm}

 Thus, if in some permutation representation of $G$, $|fix(H)|$ is too big, then $H$ is not pronormal in $G$. Therefore, it is interesting to consider pronormality of the subgroups of a group $G$ containing a subgroup $S$ which is pronormal in $G$. In particular, it is interesting to consider pronormality of overgroups of Sylow subgroups. Note that the subgroups of odd index in a finite group $G$ are exactly overgroups of Sylow $2$-subgroups of $G$.

 First of all, we will concentrate on the question of pronormality of subgroups of odd index in non-abelian simple groups. In 2012, E.~Vdovin and the third author \cite{Vdovin_Revin} formulated the following conjecture.

 \medskip

 {\bf Conjecture 1.} \emph{The subgroups of odd index are pronormal in all simple groups.}

\medskip

In this paper, we disprove Conjecture~1 and discuss a recent progress in the classification of non-abelian simple groups in which the subgroups of odd index are pronormal.

In Section~\ref{nonsimple}, we discuss the question of pronormality of subgroups of odd index in a non-simple group.
Note that an answer to this question for a group $G$ is very weakly connected with answers to this question for subgroups of $G$, even for normal subgroups of $G$. For example, it could be proved that there exist infinitely many non-abelian simple groups in which the subgroups of odd index are pronormal while a direct product of any two of them contains a non-pronormal subgroup of odd index (see Example \ref{NonPronInProd}).

\smallskip

\section{Terminology and Notation}

Our terminology and notation are mostly standard and can be found in \cite{ATLAS,Kleidman_Liebeck}.

For a group $G$ and a subset $\pi$ of the set of all primes, $O_\pi(G)$ and $Z(G)$ denote the $\pi$-radical (the largest normal $\pi$-subgroup) and the center of $G$, respectively.  As usual, $\pi'$ stands for the set of those primes that do not lie in $\pi$. If $n$ is a positive integer, then $n_\pi$ is the largest divisor of $n$ whose all prime divisors lie in $\pi$. Also, for a group $G$, it is common to write $O(G)$ instead of $O_{2'}(G)$. The set of Sylow $p$-subgroups of $G$ is denoted by $Syl_p(G)$. The socle of $G$ is denoted by $Soc(G)$. Recall that $G$ is almost simple if $Soc(G)$ is a non-abelian simple group.

We use the following notation for non-abelian simple groups: $Alt(n)$ for alternating groups;
$PSL_n^\varepsilon(q)$, where $\varepsilon=+$ for $PSL_n(q)$ and $\varepsilon=-$ for $PSU_n(q)$; $E_6^\varepsilon(q)$, where $\varepsilon=+$ for $E_6(q)$ and $\varepsilon=-$ for $^2E_6(q)$.

Let $m$ and $n$ be non-negative integers with $m=\sum_{i=0}^\infty a_i\cdot 2^i$ and $n=\sum_{i=0}^\infty b_i\cdot 2^i$, where $a_i,b_i\in\{0,1\}$. We write $m\preceq n$ if $a_i\leq b_i$ for every $i$ and $m\prec n$ if, in addition, $m \not = n$.

\section{Verification of Conjecture~1 for many families of non-abelian simple groups}

An important role in the verification of Conjecture~1 is played by the following easy assertion (see \cite[Lemma~5]{Vdovin_Revin}), which is a consequence of Theorem~\ref{Halls1Thm}.

\begin{lemma}\label{NormSyl} Suppose that $G$ is a group and $H \le G$. Assume also that $H$ contains a Sylow subgroup $S$ of $G$. Then the following statements are equivalent{\rm:}

$(1)$ $H$ is pronormal in $G${\rm;}

$(2)$ The subgroups $H$ and $H^g$ are conjugate in $\langle H, H^g \rangle$ for every $g \in N_G(S)$.

\end{lemma}

Note that Sylow $2$-subgroups in non-abelian simple groups are usually self-normalized, and all the exceptions were described by the first author in \cite{Kondrat'ev}. He proved the following proposition (see \cite[Corollary of Theorems~1-3]{Kondrat'ev}).

\begin{prop}\label{N_G(S)=S} Let $G$ be a non-abelian simple group and let $S \in Syl_2(G)$. Then $N_G(S)=S$ excluding the following cases:

$(1)$ $G \cong J_2, J_3, Suz$ or $F_5$ and $|N_G(S):S|=3${\rm;}

$(2)$ $G\cong {^2}G_2(3^{2n+1}) $ or $J_1$ and $N_G(S)\cong (C_2)^3\rtimes (C_7\rtimes C_3)<Hol((C_2)^3)${\rm;}

$(3)$ $G$ is a group of Lie type over a field of characteristic $2$ and $N_G(S)$ is a Borel subgroup of $G${\rm;}

$(4)$ $G \cong PSL_2(q)$, where $3<q\equiv \pm 3 \pmod 8$ and $N_G(S) \cong Alt(4)${\rm;}

$(5)$  $G\cong PSp_{2n}(q)$, where $n\geq 2$, $q\equiv\pm 3\pmod 8$, $n=2^{s_1}+\dots+2^{s_t}$ for \\ $s_1>\dots>s_t\geq 0$ and $N_G(S)/S$ is the elementary abelian group of order $3^t${\rm;}

$(6)$  $G \cong E^{\eta}_6(q)$ where $\eta \in \{+,-\}$, $q$ is odd, and $|N_G(S):S|=(q-\eta 1)_{2'}/(q-\eta 1,3) \not = 1;$

$(7)$ $G\cong PSL_n^\eta(q)$, where $n\geq 3$, $\eta\in \{+,-\}$, $q$ is odd, $n=2^{s_1}+\dots+2^{s_t}$ for \\ $s_1>\dots>s_t>0$, and $N_G(S)\cong S\times C_1\times \dots\times C_{t-1}>S,$ where  $C_1,\dots C_{t-2}$, and $C_{t-1}$ are cyclic subgroup of orders
${(q-\eta 1)}_{2'},\dots$, ${(q-\eta 1)}_{2'}$, and ${(q-\eta 1)}_{2'}/{(q-\eta 1,n)} _{2'}$, respectively.

\end{prop}

Thus, using Lemma~\ref{NormSyl} for the verification of Conjecture~1, it is sufficient to consider simple groups from items $(1)$--$(7)$ of Proposition~\ref{N_G(S)=S}.

In \cite{Kond_Mas_Rev1} we considered simple groups from items $(1)$--$(4)$ and, with taking into account of simple groups with self-normalized Sylow $2$-subgroups, we have proved the following theorem.

\begin{thm}\label{Array} The subgroups of odd index are pronormal in the following simple groups{\rm:}
$Alt(n)$, where $n\ge 5${\rm;}
sporadic groups{\rm;}
groups of Lie type over fields of characteristic $2${\rm;}
$PSL_{2^n}^\varepsilon(q)$, where $\varepsilon \in \{+,-\}${\rm;}
$PSp_{2n}(q)$, where  $q\not\equiv\pm3 \pmod 8${\rm;}
orthogonal groups{\rm;}
exceptional groups of Lie type not isomorphic to $E^\varepsilon_6(q)$, where $\varepsilon \in \{+,-\}$.

\end{thm}

\section{Maximal subgroups of odd index in simple groups}

In this paper, the complete classification of maximal subgroups of odd index in simple classical groups is a crucial tool in view of the following evident lemma.

\begin{lemma}\label{Overgroup} Suppose that $H$ and $M$ are subgroups of a group  $G$ and $H \le M$. Then

$(1)$ if $H$ is pronormal in $G$, then $H$ is pronormal in  $M${\rm;}

$(2)$ if $S \le H$ for some Sylow subgroup $S$ of $G$, $N_G(S) \le M$, and $H$ is pronormal in $M$, then $H$ is pronormal in $G$.

\end{lemma}

M.~Liebeck and J.~Saxl \cite{LiSa} and, independently, W.~Kantor \cite{Kantor} proposed a classification of primitive permutation groups of odd degree. It is considered to be one of remarkable results in the theory of finite permutation groups.
In particular, both papers \cite{LiSa} and \cite{Kantor} contain lists of subgroups of simple groups that can turn out to be maximal subgroups of odd index. However, in the cases of alternating groups and of classical groups over fields of odd characteristics, neither in \cite{LiSa} nor in \cite{Kantor} it was described  which of the specified subgroups are precisely maximal subgroups of odd index. Thus, the problem of the complete classification of maximal subgroups of odd index in simple groups was remained open.

The classification was finished by the second author in \cite{Maslova1, Maslova2}. In \cite{Maslova1}, the author referred to results obtained by P.~Kleidman \cite{Kleidman} and by P.~Kleidman and M.~Liebeck \cite{Kleidman_Liebeck}. However, there are a number of  inaccuracies in Kleidman's PhD thesis \cite{Kleidman}. These inaccuracies have been corrected in \cite{BrDougHolt}. Due to uncovered circumstances, in \cite{Maslova_revision} the second author revised the main result of \cite{Maslova1} taking into account the results of \cite{BrDougHolt}. In particular, changes in statements of items $(6)$, $(10)$, and $(21)$ of \cite[Theorem~1]{Maslova1} were made.

\section{Some tools for disproving Conjecture~1}

In \cite{Kond_Mas_Rev2}, we disproved Conjecture~1. The aim of this section is to demonstrate some tools for this.

 A consequence of well-known Schur--Zassenhaus theorem (see, for example, \cite[Theorems 3.8 and 3.12]{Isaacs}) is the following proposition \cite[Ch.~4, Lemma~4.28]{Isaacs}.

\begin{prop}\label{Supl} If $V$ is a normal subgroup of a group $G$ and $H$ is a subgroup of $G$ such that $(|H|,|V|)=1$, then, for any $H$-invariant subgroup $U$ of $V$, the equality  $U=C_U(H)[H,U]$ holds.

\end{prop}

We proved that the following more general statement \cite[Proposition~2]{Kond_Mas_Rev2} holds.

\begin{prop}\label{SuffSupl} If $V$ is a normal subgroup of a group $G$ and $H$ is a pronormal subgroup of $G$, then, for any $H$-invariant subgroup $U$ of $V$, the equality $U=N_U(H)[H,U]$ holds.

\end{prop}

It is easy to see that, in the case when the subgroups $H$ and $V$ from Proposition \ref{SuffSupl} have trivial intersection, the equality $N_U(H)=C_U(H)$ holds for any $H$-invariant subgroup $U$ of $V$. Therefore, Proposition \ref{Supl} is a special case of Proposition \ref{SuffSupl}.

We showed that the statement converse to Proposition \ref{SuffSupl} holds when the group $V$ is abelian and $G=HV$ (i.~e., $H$ is a supplement to the subgroup $V$ in $G$). We have proved the following theorem (see \cite[Theorem~1]{Kond_Mas_Rev2}).

\begin{thm}\label{SuffCriteria} Let $H$ and $V$ be subgroups of a group $G$ such that $V$ is an abelian normal subgroup of $G$ and $G=HV$. Then the following statements are equivalent:

$(1)$ the subgroup $H$ is pronormal in G;

$(2)$ $U=N_U(H)[H,U]$ for any $H$-invariant subgroup $U$ of $V$.

\end{thm}

With using Theorem \ref{SuffCriteria}, we have proved the following proposition (see \cite[Corollary of Theorem~1]{Kond_Mas_Rev2}).

\begin{prop}\label{CpwrSn} Let $G=A\wr Sym(n)=HV$ be the natural permutational wreath product of an abelian group $A$ and the symmetric group $H=Sym(n)$, where $V$ denotes the base of the wreath product. Then the following statements are equivalent:

$(1)$ the subgroup $H$ is pronormal in G;

$(2)$ $(|A|,n)=1$.

\end{prop}

\section{A series of examples disproving Conjecture~1}

The aim of this section is to construct a series of examples disproving Conjecture~1.
We prove the following theorem (see \cite[Theorem~2]{Kond_Mas_Rev2}).

\begin{thm}\label{Counterexample} The simple group $PSp_{6n}(q)$ for any $q \equiv \pm 3 \pmod 8$ contains a non-pronormal subgroup of odd index.

\end{thm}

{\bf Sketch of proof.} Let $q\equiv \pm 3\pmod 8$ be a prime power and $n$ be a positive integer.
It is well known that a Sylow $2$-subgroup $S$ of the group $T=Sp_2(q)=SL_2(q)$ is isomorphic to $Q_8$, and $N_T(S) \cong SL_2(3)\cong Q_8:3$.
\smallskip
We have the following chain of embeddings: $$Q_8\wr Sym(3n) \le L=Sp_2(3)\wr Sym(3n) \le Sp_2(q)\wr Sym(3n) \le  G=Sp_{6n}(q).$$

It could be proved by direct calculations that $|G:L|$ is odd.
It is easy to see that $L/O_2(L) \cong C_3 \wr Sym(3n)$. In view of Proposition~\ref{CpwrSn}, the group $L/O_2(L)$ contains a non-pronormal subgroup $T$ of odd index. Let $H$ be the preimage of $T$ in $L$. Then $H$ is a non-pronormal subgroup of odd index in $L$ in view of the following proposition.

\begin{prop}{\rm {\bf (see \cite[Lemma~3]{Mazurov} and \cite [Chapter~I, Proposition (6.4)]{Doerk_Hawks})}}\label{Quot} Suppose that $H$ is a subgroup and $N$ is a normal subgroup of a group $G$. Let $\bar{\empty} : G \rightarrow G/N$
be the natural epimorphism. The following statements hold{\rm:}

$(1)$ if $H$ is pronormal in $G$, then $\overline{H}$ is pronormal in $\overline{G}${\rm;}

$(2)$ $H$ is pronormal in $G$ if and only if $\overline{HN}$ is pronormal in $\overline{G}$ and $H$ is pronormal in $N_G(HN)${\rm;}

$(3)$ if $N \le H$ and $\overline{H}$ is pronormal in $\overline{G}$, then $H$ is pronormal in $G$.
In particular, a subgroup $H$ of odd index is pronormal in $G$ if and only if $H/O_2(G)$ is pronormal in $G/O_2(G)$.

\end{prop}

In view of Lemma~\ref{Overgroup}, $H$ is a non-pronormal subgroup of $G$, and it is easy to see that $|G:H|$ is odd. Note that $O_2(G) \le H$. Thus, $H/O_2(G)$ is a non-pronormal subgroup of odd index in $G/O_2(G) \cong PSp_{6n}(q)$ in view of Proposition~\ref{Quot}.

\section{Classification problem}

The following problem naturally arises.

\medskip

{\bf Problem 1.} \emph{Classify non-abelian simple  groups in which the subgroups of odd index are pronormal.}

\medskip

To solve Problem~1 it remains to consider the following simple groups:

$(1)$ $PSp_{2n}(q)$, where  $q \equiv\pm3 \pmod 8$ and $3$ does not divide $n${\rm;}

$(2)$ $E_6^\varepsilon(q)$, where $\varepsilon \in \{+,-\}$ and $q$ is odd{\rm;}

$(3)$ $PSL_{n}^\varepsilon(q)$, where $\varepsilon \in \{+,-\}$, $q$ is odd, and $n \not =2^w$.

\smallskip

In this paper, we consider in some details a solution of Problem 1 for symplectic groups,
briefly discuss a solution of Problem~1 for groups $E_6^\varepsilon(q)$, and formulate a conjecture for groups $PSL_n^\varepsilon(q)$.

\section{Simple symplectic groups containing non-pronormal subgroups of odd index}

In fact, Theorem \ref{Counterexample} permits us to investigate Problem~1 within a much wider family of symplectic groups.

Let $G=Sp_{2n}(q)$, where $q$ is odd, and $V$ be the natural $2n$-dimensional vector space over the field $F_q$ with a non-degenerate skew-symmetric bilinear form associated with $G$. Let us look to the list of maximal subgroups of odd index of $G$ (this list can be found in \cite{Maslova_revision}).

\begin{prop}\label{MaxSympl} Maximal subgroups of odd index in $Sp_{2n}(q)=Sp(V)$, where $n>1$ and $q$ is odd, are the following{\rm:}

 $(1)$ $Sp_{2n}(q_0)$, where $q=q_0^r$ and $r$ is an odd prime{\rm;}

 $(2)$ $Sp_{2m}(q)\times Sp_{2(n-m)}(q)$, where  $m\prec n${\rm;}

 $(3)$ $Sp_{2m}(q)\wr Sym(t)$, where  $n=mt$ and $m=2^k${\rm;}

 $(4)$ $2^{1+4}_+.Alt(5)$, where $n=2$ and $q\equiv \pm 3\pmod 8$ is a prime.

\end{prop}

Note that if $q$ is odd, then $|Z(Sp_{2n}(q))|=2$. Thus, in view of Proposition~\ref{Quot}, the subgroups of odd index are pronormal in $Sp_{2n}(q)$ if and only if the subgroups of odd index are pronormal in $PSp_{2n}(q)$, i.~e., there is no difference between investigation of Problem~1 for $Sp_{2n}(q)$ or for $PSp_{2n}(q)$.

Let $q \equiv \pm 3 \pmod8$ and $n$ is not of the form $2^w$ or $2^w(2^{2k}+ 1)$. Then the $2$-adic representation $n=\sum_{i=0}^\infty s_i\cdot 2^i$ either has two $1$s in positions $s_1$ and $s_2$ of different parity, or three $1$s in positions $s_1$, $s_2$, and $s_3$ of the same parity. Define $m=2^{s_1}+2^{s_2}$ or $m=2^{s_1}+2^{s_2}+2^{s_3}$, respectively. It is easy to see that $m\prec n$ and $3$ divides $m$.

Let $M$ be the stabilizer in $G$ of a non-degenerate $2m$-dimensional subspace of $V$. It is easy to see that $M=M_1\times M_2$, where $M_1 \cong Sp_{2m}(q)$ and $M_2\cong Sp_{2(n-m)}(q)$. Note that the index $|G:M|$ is odd by Proposition \ref{MaxSympl}. Thus, if $H$ is a subgroup of odd index in $M$, then $H$ is a subgroup of odd index in $G$ too.  Since $3$ divides $m$, it follows by Theorem \ref{Counterexample} that $M_1/Z(M_1)=M_1/O_2(M_1)$ has a non-pronormal subgroup $H_1/O_2(M_1)$ of odd index. Then $H_1$ is a non-pronormal subgroup of odd index in $M_1$ in view of Proposition~\ref{Quot}. Now it is easy to see that therefore, $H_1 \times M_2$ is a non-pronormal subgroup of odd index in $M_1 \times M_2$. So, $H_1 \times M_2$ is a non-pronormal subgroup of odd index in $G$.
Thus, we proved the following theorem (see \cite[Theorem~1]{Kond_Mas_Rev3}).

\begin{thm}\label{SimplNonPron} Let $G=PSp_{2n}(q)$, where $q \equiv \pm 3 \pmod 8$ and $n$ is not of the form $2^w$ or $2^w(2^{2k}+1)$. Then $G$ has a non-pronormal subgroup of odd index.

\end{thm}

Now in view of Theorems \ref{Array} and \ref{SimplNonPron}, to finish a classification of simple symplectic groups in which the subgroups of odd index are pronormal, it remains to consider groups $PSp_{2n}(q)$, where $q \equiv \pm 3 \pmod 8$ and $n$ is of the form $2^w$ or $2^w(2^{2k}+1)$. In the next section we will prove that the subgroups of odd index are pronormal in these groups.

\section{Classification of simple symplectic groups in which the subgroups of odd index are pronormal}

In this section we prove the following theorem, whose proof was recently finished by the authors  (see \cite[Theorem~2]{Kond_Mas_Rev3} for the case when $n$ is of the form $2^w$ and \cite{Kond_Mas_Rev4} for the case when $n$ is of the form $2^w(2^{2k}+1)$).

\begin{thm}\label{PronormalSympl} Let $G=PSp_{2n}(q)$. Then each subgroup of odd index is pronormal in $G$ if and only if one of the following statements holds:

$(1)$ $q \not \equiv \pm 3 \pmod 8${\rm;}

$(2)$ $n$ is of the form $2^w$ or $2^w(2^{2k}+1)$.

\end{thm}

{\bf Sketch of proof.} Let $G=Sp_{2n}(q)$, where $q \equiv \pm 3 \pmod 8$ and $n$ is of the form $2^w \ge 2$.
Suppose that the claim of the theorem is false, and let $q$ be the smallest prime power congruent to $\pm 3$ modulo $8$ such that $G$ has a non-pronormal subgroup $H$ of odd index. Let $S \le H$ be a Sylow $2$-subgroup of $G$. Note that $|N_G(S)/S|=3$ in view of Proposition \ref{N_G(S)=S}.

Pick $g \in N_G(S)$. Without loss of generality we can assume that $|g|=3$. Put $K=\langle H, H^g\rangle$, and we can suppose that $K < G$. Then there exists a maximal subgroup $M$ of $G$ such that $K \le M$.

Using Proposition \ref{MaxSympl}, we conclude that one of the following cases arises.

1. $M \cong Sp_{2n}(q_0)$, where $q=q_0^r$ and $r$ is an odd prime. In this case it could be proved that $N_G(S)\le M$, and we use inductive reasonings and Lemma \ref{Overgroup} to prove that $H$ is pronormal in $G$.

2. $M \cong Sp_{2m}(q)\wr Sym(t)$, where  $n=mt$ and $m=2^k$, and $M$ is choosen of such type so that $m$ is as small as possible.
In this case we can prove that $H/O_2(M)$ is pronormal in $M/O_2(M) \cong PSp_{2m}(q)\wr Sym(t)$ with using two following propositions and additional inductive reasonings.

\begin{prop}{\rm {\bf {\rm {\bf (see \cite[Lemma~15]{Kond_Mas_Rev3}}} and \cite[Lemma~9]{Guo_Mas_Rev})}}\label{DirProd}
Let $Q$ be a subgroup of odd index in a group $L=L_1\times L_2 \times \ldots \times L_n$, where $L_i$ are groups, and let $\pi_i$ be the projection from $L$ to $L_i$. If there is $i$ such that $L_i$ is almost simple, $L_i/ Soc(L_i)$
is a $2$-group, and $Q^{\pi_i}=L_i$, then $L_i \le Q$.
\end{prop}

\begin{prop}{\rm {\bf {\rm {\bf (see \cite[Lemma~17]{Kond_Mas_Rev3}}})}}\label{PSpwrSn} Let $G=A\wr Sym(n)=LH$ be the natural permutation wreath product of a non-abelian simple group $A$ and $H=Sym(n)$, where $n=2^w$, $L=L_1\times \ldots \times L_n$, and for each $i \in \{1, \ldots, n\}$, $L_i \cong A$ and $\pi_i: L \rightarrow L_i$ is the projection from $L$ to $L_i$. If $K$ is a subgroup of odd index in $G$, $K_0=K \cap L$, and $M_1 \le L_1$ such that $N_{L_1}(K_0^{\pi_1})\le M_1$, then $K \le U \cong M_1 \wr  Sym(n) \le G$.

\end{prop}

\smallskip

Moreover, it could be proved that in this case $N_G(S)\le M$. Now, in view of  Proposition~\ref{Quot} and Lemma~\ref{Overgroup}, we conclude that $H$ is pronormal in $G$.

3. $M \cong 2^{1+4}_+.Alt(5)$, where $n=4$ and $q$ is a prime. In this case it could be easy proved that $N_G(S)\le M$. Now we conclude that $H$ is pronormal in $G$ in view of Proposition~\ref{Quot}, Lemma~\ref{Overgroup}, and Theorem \ref{Array}.
\medskip

In the case when $n$ is of the form $2^w(2^{2k}+1)$, the scheme of the proof is similar, however we can not use Proposition \ref{PSpwrSn}. The proof of Proposition \ref{PSpwrSn} is based on a fact that $K$ contains a regular subgroup of some conjugate of $H$. And the fact could be false if $n$ is not of the form $2^w$.

To investigate this case the second and the third authors in a joint work with W.~Guo \cite{Guo_Mas_Rev} obtained the following useful criterion of the pronormality of subgroups of odd index in extensions of groups (see \cite[Theorem~1]{Guo_Mas_Rev}).

\begin{thm}\label{CritExten} Let $G$ be a group, $A\unlhd G$,
the subgroups of odd index are pronormal in $A$, and Sylow $2$-subgroups of $G/A$ are self-normalized. Let $T$ be a Sylow $2$-subgroup of $A$. Then the following statements are equivalent{\rm:}

  $(1)$ The subgroups of odd index are pronormal in $G${\rm;}

  $(2)$ The subgroups of odd index are pronormal in  $N_G(T)/T$.

\end{thm}

However, there are difficulties in a direct application of Theorem \ref{CritExten} to the maximal subgroups of the group $Sp_{2n}(q)$, where $n$ is of the form $2^w(2^{2k}+1)$, in view of the following fact: if the subgroups of odd index are pronormal in both groups $G_1$ and $G_2$, it does not imply that the subgroups of odd index are pronormal in the group $G_1 \times G_2$, even if $G_1$ and $G_2$ are non-abelian simple. It is clear from the following example (see \cite[Proposition~1]{Guo_Mas_Rev}).

\begin{ex}\label{NonPronInProd} Consider Frobenius groups $H_i = L_i\rtimes K_i \cong C_7 \rtimes C_3$ for $i \in \{1,2\}$ and $H = H_1 \times H_2$. Note that any proper subgroup of $H_i$ is its Sylow subgroup, and hence is pronormal in $H_i$.  Let $L=O_7(H)=L_1\times L_2 \cong C_7 \times C_7$ and $D=\{(x,x) \mid x \in C_7\}.$

There exists $k_1 \in K_1 \times \{1\}$ such that $D^{k_1}=\{(2x,x) \mid x \in C_7\}.$ Hence $\langle D, D^{k_1}\rangle=L$ is abelian and $D$ is a non-pronormal subgroup {\rm(}of odd index{\rm)} in $H$.

Let $G_1, G_2 \in \{J_1\} \cup \{{^2}G_2(3^{2m+1})\mid m \ge 1\}$ and $S_i \in Syl_2(G_i)$. In view of Theorem \ref{Array}, the subgroups of odd index are pronormal in $G_i$ for $i \in \{1,2\}$. And in view of Proposition~\ref{N_G(S)=S}, $N_{G_i}(S_i)/S_i$  for $i \in \{1,2\}$ is isomorphic to the Frobenius group $C_7 \rtimes C_3$. Using previous reasonings it is easy to construct a non-pronormal subgroup of odd index in $G=G_1 \times G_2$.

\end{ex}


\smallskip

Recently basing on Theorem \ref{CritExten}, the second and the third authors, and W.~Guo \cite{Guo_Mas_Rev} obtained the following pronormality criterion for subgroups of odd index in groups of the type $\prod_{i=1}^t (A\wr Sym(n_i))$, where $A$ is an abelian group and all the wreath products are natural permutation (see \cite[Theorem~2]{Guo_Mas_Rev}).

\begin{thm}\label{AwrSn} Let $A$ be an abelian group and $G=\prod_{i=1}^t (A\wr Sym(n_i))$, where all the wreath products are natural permutation. Then the subgroups of odd index are pronormal in $G$ if and only if for any positive integer $m$, the inequality $m \prec  n_i$ for some $i$ implies that $g.c.d.(|A|,m)$ is a power of $2$.

\end{thm}

Moreover, in \cite[Theorem~3]{Guo_Mas_Rev} the following theorem was proved.

\begin{thm}\label{ProdSympl} Let $G =\prod_{i=1}^t G_i$, where for any $i \in \{1, . . . , t\}$, $G_i \cong Sp_{2n_i}(q_i)$, $q_i$ is odd, and $n_i$ is a power of $2$. Then the subgroups of odd index are pronormal in $G$.

\end{thm}

Theorems \ref{CritExten}, \ref{AwrSn}, and \ref{ProdSympl} became the main tools in the proof of Theorem~\ref{PronormalSympl} which was finished by the authors in \cite{Kond_Mas_Rev4}.

\section{Summary and further research on finite simple groups in which the subgroups of odd index are pronormal}

Let $G$ be a non-abelian simple  group, $S \in Syl_2(G)$, and $C=O(C_G(S))$.
In \cite[Theorem~7]{KondMaz} it was proved that $C \not =1$ if and only if one of the following statements holds:

$(1)$ $G \cong E^{\eta}_6(q)$, where $\eta \in \{+,-\}$, $q$ is odd, and $C$ is a cyclic group of order $(q-\eta 1)_{2'}/(q-\eta 1,3) \not = 1${\rm;}

$(2)$ $G\cong PSL_n^\eta(q)$, where $n\geq 3$, $\eta\in \{+,-\}$, $q$ is odd, $n=2^{s_1}+\dots+2^{s_t}$ for \\ $s_1>\dots>s_t>0$ and $t>1$, and $C\cong C_1\times \dots\times C_{t-1} \not =1,$ where  $C_1,\dots C_{t-2}$, and $C_{t-1}$ are cyclic subgroup of orders
${(q-\eta 1)}_{2'},\dots$, ${(q-\eta 1)}_{2'}$, and ${(q-\eta 1)}_{2'}/{(q-\eta 1,n)} _{2'}$, respectively.

\smallskip

Thus, with taking into account of Theorems \ref{Array} and \ref{PronormalSympl}, we obtain the following theorem.

\begin{thm} Let $G$ be a non-abelian simple group, $S \in Syl_2(G)$, and $C=O(C_G(S))$. If $C = 1$, then exactly one of the following statements holds{\rm:}

$(1)$ The subgroups of odd index are pronormal in $G${\rm;}

$(2)$ $G\cong PSp_{2n}(q)$, where $q \equiv \pm 3 \pmod 8$ and $n$ is not of the form $2^w$ or $2^w(2^{2k}+1)$.

\end{thm}

Moreover, recently we have obtained a solution of Problem 1 for simple exceptional groups $E_6^\varepsilon(q)$. We have proved the following theorem.

\begin{thm}\label{E_6} Let $G=E_6^\varepsilon(q)$, where $q=p^k$, $p$ is a prime, and $\varepsilon \in \{+,-\}$. Then the subgroups of odd index are pronormal in $G$ if and only if $18$ does not divide $q-\varepsilon1$ and if $p$ is odd and $\varepsilon=+$, then $k$ is a power of $2$.
\end{thm}

The scheme of our proof of Theorem \ref{E_6} is similar as for symplectic groups and is based on the classification of maximal subgroups of odd index in simple exceptional groups of Lie type obtained in \cite{LiSa,Kantor}. Moreover, we use some results on subgroup structure of $G=E_6^\varepsilon(q)$ obtained in \cite{KondMaz}. The most difficult case here is when a possibly non-pronormal subgroup $H$ of odd index and its conjugate $H^g$ for some $g\in G$ are both contained in a parabolic maximal subgroup of $G$. In this case Theorem~\ref{CritExten} is an useful tool.

\medskip

Thus, to solve Problem 1 it remains to consider the groups $PSL_{n}^\varepsilon(q)$, where $\varepsilon \in \{+,-\}$, $q$ is odd, and $n \not =2^w$. We have the following conjecture.

\medskip

{\bf Conjecture 2.} \emph{Let $G=PSL_n^\varepsilon(q)$, where $q$ is odd and $\varepsilon \in \{+,-\}$.
 The subgroups of odd index are pronormal in $G$ if and only if for any positive integer $m$, the inequality $m \prec n$ implies that $g.c.d.(m,q^{1+\varepsilon1}(q-\varepsilon1))$ is a power of $2$.}

\medskip

\section{Question of the pronormality of subgroups of odd index in non-simple groups}\label{nonsimple}

The Frattini Argument (see Proposition~\ref{Frattini}) and Proposition~\ref{Quot} are convenient tools to reduce General Problem to groups of smaller order.
Assume that $G$ is not simple and $A$ is a minimal non-trivial normal subgroup of $G$.
Then $A$ is a direct product of pairwise isomorphic simple groups, and one of the following cases arises:

$(1)$ $A \le H$ and, in view of Proposition~\ref{Quot}, $H$ is pronormal in $G$ if and only if $H/A$ is pronormal in $G/A$.
Note that $|G/A|<|G|$.

$(2)$ $H \le A$ and, in view of Proposition~\ref{Frattini}, $H$ is pronormal in $G$ if and only if $H$ is pronormal in $A$ and $G=AN_G(H)$. Note that $|A|<|G|$. Thus, the question of pronormality of subgroups in direct products of simple groups is of interest.

$(3)$ Let $H \not\le A$ and $A \not\le H$, and let $N=N_G(HA)$. In view of Propositions~\ref{Quot}~and~\ref{Frattini}, $H$ is pronormal in $G$ if and only if $HA/A$ is pronormal in $G/A$, $N=AN_N(H)$, and $H$ is pronormal in $HA$. Therefore, with using inductive reasonings, we can reduce this case to the subcase when $G=HA$.

\smallskip

Suppose that $G=HA$, $A$ is a minimal non-trivial normal subgroup of $G$, $A \not \le H$, and $|G:H|$ is odd.

If $|A|$ is odd, then $A$ is abelian and $H$ is pronormal in $G$. Indeed, if $U$ is an $H$-invariant subgroup of $A$, then either $U$ is trivial or $U=A$. Therefore, $A \cap H$ is trivial, and it is easy to see that $U=C_U(H)[H,U]=N_U(H)[H,U]$ for every $H$-invariant subgroup $U$ of $A$. Thus, in view of Theorem~\ref{SuffCriteria}, $H$ is pronormal in $G$.

If $A$ is a $2$-group, then $H$ is pronormal in $G$ in view of Proposition~\ref{Quot}.

Suppose that $A$ is a direct product of pairwise isomorphic non-abelian simple groups. If the subgroups of odd index are pronormal in $A$, then
we can use the following criterion of the pronormality of subgroups of odd index in extensions of groups, which is a generalization of Theorem~\ref{CritExten}.

\begin{thm}\label{CritExtenExt} Let $G$ be a group, $A\unlhd G$, and the subgroups of odd index are pronormal both in $A$ and in $G/A$. Let $T$ be a Sylow $2$-subgroup of $A$ and let $\bar{\empty} : G \rightarrow G/A$ be the natural epimorphism.

$(1)$ Assume that $H \le G$, $S\le H$ for some $S\in Syl_2(G)$, and $T = A\cap S$. Let $Y=N_{A}(H\cap A)$ and $Z=N_{H\cap A}(T)$.
Then $H$ is pronormal in $G$ if and only if $N_H(T)/Z$  is pronormal in $(N_H(T) N_Y(T) )/Z$ and $\overline{N_G(H)} =N_{\overline{G}}(\overline{H})$.

$(2)$ The subgroups of odd index are pronormal in $G$ if and only if the subgroups of odd index are pronormal in  $N_G(T)/T$ and
for every subgroup $H \le G$ of odd index, we have $\overline{N_G(H)} =N_{\overline{G}}(\overline{H})$.
\end{thm}

Note that proof of Theorem~\ref{CritExtenExt} follows from proofs of Theorems~1 and~4 in \cite{Guo_Mas_Rev}.

\medskip

In view of Theorem~\ref{CritExtenExt}, the following problem is of interest.

\medskip

{\bf Problem 2.} \emph{Describe direct products of non-abelian simple groups in which the subgroups of odd index are pronormal.}

\medskip

Note that Problem~2 is not equivalent to Problem~1 (see Example~\ref{NonPronInProd}). However, in some cases the pronormality of subgroups of odd index in a direct product of non-abelian simple groups is equivalent to the pronormality of subgroups of odd index in each factor. As an example, recently we  have proved the following theorem.

\begin{thm}\label{DirProd} Let $G =\prod_{i=1}^t G_i$, where $G_i \cong PSp_{2n_i}(q_i)$ and $q_i$ is odd for each $i \in \{1, . . . , t\}$. Then the subgroups of odd index are pronormal in $G$ if and only if the subgroups of odd index are pronormal in $G_i$ for each $i$.
\end{thm}

\bigskip




\end{document}